\documentclass[12pt]{article}
\usepackage{amssymb, amsmath, url, graphicx}

\def\3{\subset }
\def\4{\subseteq }
\def\<{\left<}
\def\>{\right>}
\def\vsp{\vspace*{1,5mm}\\ }

\def\bit{\begin{itemize}}
\def\eit{\end{itemize}}
\def\3{\subset }
\def\4{\subseteq }
\def\ov{\overline}

\def\0{\leqno}

\def\barr{\begin{array}}
\def\earr{\end{array}}
\def\dd{\displaystyle}

\def\Z{{\rlap{$\kern2pt{\rm Z}$}{\rm Z}\,}}
\def\bld#1#2{{\buildrel{#1}\over{#2}}}
\def\st#1#2{{\mathrel{\mathop{#2}\limits_{#1}}{}\!}}
\def\stb#1#2#3{{\st{{#1}}{\bld{{#2}}{#3}}{}\!}}
\def\xmare#1#2{\stb{#1}{#2}{\mbox{\Huge$\times$}}}

\def\calg{{\cal G}}

\def\frax{\dd\frac}
\def\io{if and only if }

\title{\bf Cyclic subgroup commutativity degrees of finite groups}
\author{Marius T\u arn\u auceanu and Mihai-Silviu Lazorec}
\date{September 2, 2016}

\begin{document}

\maketitle

\begin{abstract}
In this paper we introduce and study the concept of cyclic
subgroup commutativity degree of a finite group $G$. This quantity
measures the probability of two random cyclic subgroups of $G$
commuting. Explicit formulas are obtained for some particular
classes of groups. A criterion for a finite group to be an Iwasawa group
is also presented.
\end{abstract}



\section{Introduction}

In the last years there has been a growing interest in the use of
probability in finite group theory. One of the most important
aspects which have been studied is the probability that two
elements of a finite group $G$ commute. It is called the {\it
commutativity degree} of $G$ and has been investigated in many
papers, such as \cite{2,3} and \cite{5}--\cite{9}. Inspired by
this concept, in \cite{16} we introduced a similar notion for the
subgroups of $G$, called the {\it subgroup commutativity degree}
of $G$. This quantity is defined by
$$\barr{lcl} sd(G)&=&\frax1{|L(G)|^2}\,\left|\{(H,K)\in L(G)^2\mid
HK=KH\}\right|=\vsp &=&\frax1{|L(G)|^2}\,\left|\{(H,K)\in
L(G)^2\mid HK\in L(G)\}\right|\earr\0(1)$$ (where $L(G)$ denotes
the subgroup lattice of $G$) and it measures the probability that
two subgroups of $G$ commute, or equivalently the probability that
the product of two subgroups of $G$ be a subgroup of $G$ (recall
also the natural generalization of $sd(G)$, namely the {\it
relative subgroup commutativity degree} of a subgroup of $G$,
introduced and studied in \cite{17}).
\bigskip

Another two probabilistic notions on $L(G)$ have been investigated
in \cite{21} and \cite{18}: the {\it normality degree} and the
{\it cyclicity degree} of $G$. They are defined by
$$ndeg(G)=\dd\frac{|N(G)|}{|L(G)|}\hspace{1mm}\mbox{ and }\hspace{1mm} cdeg(G)=\dd\frac{|L_1(G)|}{|L(G)|}\,,$$where $N(G)$
and $L_1(G)$ denote the normal subgroup lattice and the poset of
cyclic subgroups of $G$, and measure the probability of a random
subgroup of $G$ to be normal or cyclic, respectively.
\bigskip

Clearly, in the definition of $sd(G)$ we may restrict to one of
the above remarkable subsets of $L(G)$. In the case of $N(G)$
nothing can be said, since normal subgroups commute with all
subgroups of $G$. By taking $L_1(G)$ instead of $L(G)$ in (1) a
new significant quantity is obtained, namely

$$\barr{lcl} csd(G)&=&\frax1{|L_1(G)|^2}\,\left|\{(H,K)\in L_1(G)^2\mid
HK=KH\}\right|=\vsp &=&\frax1{|L_1(G)|^2}\,\left|\{(H,K)\in
L_1(G)^2\mid HK\in L(G)\}\right|.\earr$$This measures the
probability that two cyclic subgroups of $G$ commute and will be
called the {\it cyclic subgroup commutativity degree} of $G$. Its
study is the purpose of the current paper.
\bigskip

The paper is organized as follows. Some basic properties and
results on cyclic subgroup commutativity degree are presented in
Section 2. Section 3 deals with cyclic subgroup commutativity
degrees for some special classes of finite groups: $P$-groups,
dihedral groups and $p$-groups possessing a cyclic maximal
subgroup. As an application, in Section 4 we give a criterion
for a finite group to be an Iwasawa group. In the final section some
further research directions and a list of open problems are indicated.
\bigskip

Most of our notation is standard and will usually not be repeated
here. Elementary notions and results on groups can be found in
\cite{4,14}. For subgroup lattice concepts we refer the reader to
\cite{13,15,20}.

\section{Basic properties of cyclic subgroup\\ commutativity degree}

Let $G$ be a finite group. First of all, we remark that the cyclic
subgroup commutativity degree $csd(G)$ satisfies the following
relation $$0<csd(G)\le 1.$$Moreover, by consequence (9) on page
202 of \cite{13}, the permutability of a subgroup $H\in L_1(G)$
with all cyclic subgroups of $G$ is equivalent with the
permutability of $H$ with all subgroups of $G$. This shows that
$$csd(G)=1\Longleftrightarrow sd(G)=1$$and therefore the finite
groups $G$ satisfying $csd(G)=1$ are in fact the Iwasawa groups,
i.e. the nilpotent modular groups (see [13, Exercise 3, p. 87]).
Notice that a complete description of these groups is given by a
well-known Iwasawa's result (see Theorem 2.4.13 of \cite{13}). In
particular, we infer that $csd(G)=1$ for all Dedekind groups $G$.
\bigskip

Given $H\in L_1(G)$, we will denote by $C_1(H)$ the set consisting
of all cyclic subgroups of $G$ commuting with $H$, that is
$$C_1(H)=\{K\in L_1(G)\mid HK=KH\}.$$Then
$$csd(G)=\frax1{|L_1(G)|^2}\dd\sum_{H\in L_1(G)}|C_1(H)|,\0(2)$$which
leads to a precise expression of $csd(G)$ for finite groups $G$
whose cyclic subgroup structure is known.

\bigskip\noindent{\bf Example 2.1.} The alternating group $A_4$ has
eight cyclic subgroups, namely: the trivial subgroup $H_1$, three
subgroups $H_i\cong\mathbb{Z}_2$, $i=2,3,4$, and four subgroups
$H_i\cong\mathbb{Z}_3$, $i=5,6,7,8$. We can easily see that
$|C_1(H_1)|=8$, $|C_1(H_i)|=4$ for $i=\ov{2,4}$, and
$|C_1(H_i)|=5$ for $i=\ov{5,8}$. Hence
$$csd(A_4)=\frax{1}{64}\left(8+3\cdot4+4\cdot5\right)=\frax{5}{8}\,.$$
\bigskip

Clearly, we have $L(H)\cup\left(N(G)\cap L_1(G)\right)\subseteq
C_1(H)$, $\forall\hspace{1mm} H\in L_1(G)$, implying that
$$csd(G)\geq\frax1{|L_1(G)|^2}\dd\sum_{H\in L_1(G)}|L(H)\cup\left(N(G)\cap
L_1(G)\right)|.$$By this inequality some lower bounds for $csd(G)$
can be inferred, namely
$$csd(G)\geq\frax1{|L_1(G)|^2}\dd\sum_{H\in L_1(G)}|N(G)\cap
L_1(G)|=\frax{|N(G)\cap L_1(G)|}{|L_1(G)|}$$and
$$csd(G)\geq\frax1{|L_1(G)|^2}\dd\sum_{H\in
L_1(G)}|L(H)|\geq\frax{2|L_1(G)|-1}{|L_1(G)|^2}\,,$$since
$|L(H)|\geq 2$\, for every non-trivial cyclic subgroup $H$ of $G$.
Another lower bound for $csd(G)$ follows by the simple remark that
for every subgroup $M$ of $G$ we have
$$\{(H,K)\in L_1(G)^2\mid HK=KH\}\supseteq\{(H,K)\in L_1(M)^2\mid
HK=KH\}.$$Thus
$$csd(G)\geq\left(\frax{|L_1(M)|}{|L_1(G)|}\right)^2 csd(M).$$In
particular, if $M$ is abelian, then $csd(M)=1$ and so
$$csd(G)\geq\left(\frax{|L_1(M)|}{|L_1(G)|}\right)^2.$$

Assume next that $G$ and $G'$ are two finite groups. If $G\cong
G'$, then $csd(G)=csd(G')$. The same thing cannot be said in the
case when $G$ and $G'$ are only lattice-isomorphic, as shows the
following elementary example.

\bigskip\noindent{\bf Example 2.2.} It is well-known that the subgroup lattices of
$G=\mathbb{Z}_3\times\mathbb{Z}_3$ and $G'=S_3$ are isomorphic. On
the other hand, we have $csd(G)=1$ because $G$ is abelian, but
$csd(G')\neq 1$ because $G'$ is not nilpotent (more precisely, we
can easily check that $csd(G')=19/25$).
\bigskip

By a direct calculation, one obtains
$$csd(S_3\times\mathbb{Z}_3)=\frax{85}{121}\ne\frax{19}{25}=csd(S_3)csd(\mathbb{Z}_3)$$and
consequently in general we don't have $csd(G\times
G')=csd(G)csd(G')$. A sufficient condition in order to this
equality holds is that $G$ and $G'$ be of coprime orders. This
remark can naturally be extended to arbitrary finite direct
products.

\bigskip\noindent{\bf Proposition 2.3.} {\it Let
$(G_i)_{i=\overline{1,k}}$ be a family of finite groups having
coprime orders. Then
$$csd(\xmare{i=1}k G_i)=\prod_{i=1}^k csd(G_i).\0(3)$$}
\bigskip

The following immediate consequence of Proposition 2.3 shows that
computing the cyclic subgroup commutativity degree of a finite
nilpotent group is reduced to finite $p$-groups.

\bigskip\noindent{\bf Corollary 2.4.} {\it If $G$ is a finite
nilpotent group and $(G_i)_{i=\ov{1,k}}$ are the Sylow subgroups
of $G$, then $$csd(G)=\prod_{i=1}^k csd(G_i).$$}

\bigskip\noindent{\bf Remark 2.5.} The condition in the hypothesis of
Proposition 2.3 is not necessary to obtain the equality (3). For
example, we have
$$csd(S_3\times\mathbb{Z}_2)=\frax{19}{25}=csd(S_3)csd(\mathbb{Z}_2),$$even
if the groups $S_3$ and $\mathbb{Z}_2$ are not of coprime orders.

\section{Cyclic subgroup commutativity degrees\\ for some classes of finite groups}

In this section we will compute explicitly the cyclic subgroup
commutativity degree of several semidirect products for which we
are able to describe the cyclic subgroup structure.
\bigskip

\bigskip
\noindent{\bf 3.1. The cyclic subgroup commutativity degree of
finite $P$-groups}
\bigskip

First of all, we recall the notion of $P$-group, according to
\cite{13}. Let $p$ be a prime, $n\geq 2$ be a cardinal number and
$G$ be a group. We say that $G$ belongs to the class $P(n,p)$ if
it is either elementary abelian of order $p^n$, or a semidirect
product of an elementary abelian normal subgroup $M$ of order
$p^{n-1}$ by a group of prime order $q\neq p$ which induces a
nontrivial power automorphism on $M$. The group $G$ is called a
$P$-$group$ if $G\in P(n,p)$ for some prime $p$ and some cardinal
number $n\geq 2$. It is well-known that the class $P(n,2)$
consists only of the elementary abelian group of order $2^n$.
Also, for $p>2$ the class $P(n,p)$ contains the elementary abelian
group of order $p^n$ and, for every prime divisor $q$ of $p-1$,
exactly one non-abelian $P$-group with elements of order $q$.
Moreover, the order of this group is $p^{n-1}q$ if $n$ is finite.
The most important property of the groups in a class $P(n,p)$ is
that they are all lattice-isomorphic (see Theorem 2.2.3 of
\cite{13}).
\bigskip

In the following, we will focus on finite non-abelian $P$-groups.
So, assume that $p>2$ and $n\in\mathbb{N}$ are fixed, and take a
divisor $q$ of $p-1$. The non-abelian group of order $p^{n-1}q$ in
the class $P(n,p)$ will be denoted by $G_{n,p}$. By Remarks 2.2.1
of \cite{13}, it is of type
$$G_{n,p}=M\langle x\rangle,$$ where $M\cong \mathbb{Z}_p^{n-1}$ (i.e. the direct product of $n-1$ copies
of $\mathbb{Z}_p$), $o(x)=q$ and there exists an integer $r$ such
that $x^{-1}yx=y^r$, for all $y\in M$. Notice that we have
$$N(G_{n,p})=L(M)\cup\{G_{n,p}\}.$$The set $L_1(G_{n,p})$ has been described in \cite{15}: it consists of the trivial subgroup
1, of the subgroups of order $p$ in $M$ and of the subgroups of
type $\langle yx\rangle$ with $y\in M$. Then
$$|L_1(G_{n,p})|=1+\frax{p^{n-1}-1}{p-1}+p^{n-1}=2+p+p^2+...+p^{n-1}.$$On the other hand, we
have $$C_1(H)=L_1(G_{n,p}), \mbox{ for all } H\leq M,$$ and
$$C_1(\langle yx\rangle)=L_1(M)\cup\{\langle yx\rangle\}, \mbox{ for
all } y\in M.$$In this way, an explicit value of $csd(G_{n,p})$ is
obtained by using (2).

\bigskip\noindent{\bf Theorem 3.1.1.} {\it The cyclic subgroup commutativity degree
of the $P$-group $G_{n,p}$ is given by the following equality:
$$csd(G_{n,p}){=}\frax{\left(2{+}p{+}p^2{+}...{+}p^{n-2}\right)\hspace{-0,5mm}\left(2{+}p{+}p^2{+}...{+}p^{n-1}\right)\hspace{-0,5mm}{+}p^{n-1}\hspace{-0,5mm}\left(3{+}p{+}p^2{+}...{+}p^{n-2}\right)}{\left(2{+}p{+}p^2{+}...{+}p^{n-1}\right)^2}\,.$$}
\smallskip

We observe that for $p=3$, $q=2$ and $n=2$ we have $G_{2,3}\cong
S_3$, and hence $csd(S_3)=\frax{19}{25}$ can be also computed by
the above formula. The following consequence of Theorem 3.1.1 is
immediate, too.

\bigskip\noindent{\bf Corollary 3.1.2.} {\it $\dd\lim_{n\to\infty}csd(G_{n,p})=\frax{2}{p}\,.$}
\bigskip

\bigskip
\noindent{\bf 3.2. The cyclic subgroup commutativity degree of
finite dihedral groups}
\bigskip

The dihedral group $D_{2m}$ $(m\ge 2)$ is the symmetry group of a
regular polygon with $m$ sides and it has the order $2m$. The most
convenient abstract description of $D_{2m}$ is obtained by using
its generators: a rotation $x$ of order $m$ and a reflection $y$
of order $2$. Under these notations, we have
$$D_{2m}=\langle x,y\mid x^m=y^2=1,\ yxy=x^{-1}\rangle.$$It is well-known that for every divisor $r$ or $m$, $D_{2m}$ possesses a subgroup isomorphic to $\mathbb{Z}_r$, namely $H^r_0=\langle x^{\frac{m}{r}}\rangle$, and $\frac{m}{r}$ subgroups isomorphic to $D_{2r}$, namely $H^r_i=\langle x^{\frac{m}{r}},x^{i-1}y\rangle,$ $i=1,2,...,\frac{m}{r}\hspace{1mm}.$ The
normal subgroups of $D_{2m}$ are
$$N(D_{2m})=\left\{\barr{lll}
L(H^m_0)\cup\{G\},& m\equiv 1 \hspace{1mm}({\rm mod}\hspace{1mm} 2)\\
\\
L(H^m_0)\cup\{G, H^{\frac{m}{2}}_1, H^{\frac{m}{2}}_2\},& m\equiv
0 \hspace{1mm}({\rm mod}\hspace{1mm} 2),\earr\right.$$while the
cyclic subgroups of $D_{2m}$ are
$$L_1(D_{2m})=L(H^m_0)\cup\{H^1_i\mid i=1,2,...,m\}.$$It follows
that $$|L_1(D_{2m})|=\tau(m)+m,$$where $\tau(m)$ denotes the
number of divisors of $m$. Clearly, we have
$$|C_1(H)|=\tau(m)+m, \mbox{ for all } H\in L(H^m_0).$$On the
other hand, it is easy to see that
$$C_1(H^1_i)=\left\{\barr{lll}
L(H^m_0)\cup\{H^1_i\},& m\equiv 1 \hspace{1mm}({\rm mod}\hspace{1mm} 2)\\
\\
L(H^m_0)\cup\{H^1_i, H^1_{i+\frac{m}{2}}\},& m\equiv 0
\hspace{1mm}({\rm mod}\hspace{1mm} 2)\earr\right.$$and therefore
$$|C_1(H^1_i)|=\left\{\barr{lll}
\tau(m)+1,& m\equiv 1 \hspace{1mm}({\rm mod}\hspace{1mm} 2)\\
\\
\tau(m)+2,& m\equiv 0 \hspace{1mm}({\rm mod}\hspace{1mm}
2),\earr\right.$$for all $i=1,2,...,m$. Then (2) leads to the
following result.

\bigskip\noindent{\bf Theorem 3.2.1.} {\it The cyclic subgroup commutativity degree
of the dihedral group $D_{2m}$ is given by the following equality:
$$csd(D_{2m})=\left\{\barr{lll}
\frax{\tau(m)(\tau(m)+m)+m(\tau(m)+1)}{(\tau(m)+m)^2}\,,& m\equiv 1 \hspace{1mm}({\rm mod}\hspace{1mm} 2)\\
\\
\frax{\tau(m)(\tau(m)+m)+m(\tau(m)+2)}{(\tau(m)+m)^2}\,,& m\equiv
0 \hspace{1mm}({\rm mod}\hspace{1mm} 2)\,.\earr\right.$$}
\smallskip

The cyclic subgroup commutativity degree of the dihedral group
$D_{2^n}$ is obtained directly from Theorem 3.2.1.

\bigskip\noindent{\bf Corollary 3.2.2.} {\it We have
$$csd(D_{2^n})=\frax{n^2+(n+1)2^n}{(n+2^{n-1})^2}$$and in particular
$$csd(D_8)=\frax{41}{49}\,.$$}
\bigskip

We are also able to compute the limit of $csd(D_{2^n})$ when
$n\to\infty$.

\bigskip\noindent{\bf Corollary 3.2.3.} {\it $\dd\lim_{n\to\infty}csd(D_{2^n})=0.$}
\bigskip
\newpage

\bigskip
\noindent{\bf 3.3. The subgroup commutativity degree of finite
$p$-groups\\ possessing a cyclic maximal subgroup}
\bigskip

Let $p$ be a prime, $n\ge 3$ be an integer and denote by $\calg$
the class consisting of all finite $p$-groups of order $p^n$
having a maximal subgroup which is cyclic. Obviously, $\calg$
contains finite abelian $p$-groups of type
$\mathbb{Z}_p\times\mathbb{Z}_{p^{n-1}}$ whose cyclic subgroup
commutativity degree is $1$, but some finite non-abelian $p$-groups
belong to $\calg$, too. They are exhaustively described in Theorem
4.1, \cite{14}, II: a non-abelian group is contained in $\calg$ \io
it is isomorphic to $M(p^n)$ when $p$ is odd, or to one of the
following groups
\begin{itemize}
\item[--] $M(2^n)\ (n\ge 4),$
\item[--] the dihedral group $D_{2^n}$,
\item[--] the generalized quaternion group
$$Q_{2^n}=\langle x,y\mid x^{2^{n-1}}=y^4=1,\ yxy^{-1}=x^{2^{n-1}-1}\rangle,$$
\item[--] the quasi-dihedral group
$$S_{2^n}=\langle x,y\mid x^{2^{n-1}}=y^2=1,\ y^{-1}xy=x^{2^{n-2}-1}\rangle\,\, (n\ge 4)$$
\end{itemize}
when $p=2$.
\bigskip

In the following the cyclic subgroup commutativity degrees of the
above $p$-groups will be determined. As we observed in Section 2,
we have $$csd(M(p^n))=1.$$Because $csd(D_{2^n})$ has been obtained
in 3.2, we need to focus only on computing $csd(Q_{2^n})$ and
$csd(S_{2^n})$.

\bigskip\noindent{\bf Theorem 3.3.1.} {\it The cyclic subgroup commutativity degree
of the generalized quaternion group $Q_{2^n}$ is
$$csd(Q_{2^n})=\frax{n^2+(n+1)2^{n-1}}{(n+2^{n-2})^2}\,.$$In particular, we have
$$csd(Q_{16})=\frax{7}{8}\,.$$}
\smallskip

\noindent{\bf Proof.} Under the above notation, it is easy to see
that $L_1(Q_{2^n})$ consists of all subgroups contained in
$\langle x\rangle$ and of all subgroups of type $\langle
x^ky\rangle$, $k=0,1,...,2^{n-2}-1$. Moreover, we have
$$|C_1(H)|=|L_1(Q_{2^n})|=n+2^{n-2},
\forall\hspace{1mm}H\leq\langle x\rangle.$$We also remark that
$$\langle x^{k_1}y\rangle\langle x^{k_2}y\rangle=\langle x^{k_2}y\rangle\langle
x^{k_1}y\rangle\Longleftrightarrow k_1=k_2 \mbox{ or }
|k_1-k_2|=2^{n-3}.$$This leads to $$|C_1(\langle
x^ky\rangle)|=n+2, \,\forall\hspace{1mm}k=0,1,...,2^{n-2}-1.$$One
obtains \bigskip
$$\hspace{-30mm}csd(Q_{2^n})=\frax{1}{(n+2^{n-2})^2}\dd\sum_{H\in L_1(Q_{2^n})}|C_1(H)|=$$
$$\hspace{27mm}=\frax{1}{(n+2^{n-2})^2}\left[\dd\sum_{^{H\in
L_1(Q_{2^n})}_{\hspace{3mm}H\leq\langle
x\rangle}}|C_1(H)|+\dd\sum_{k=0}^{2^{n-2}-1}|C_1(\langle
x^ky\rangle)|\right]=$$
$$\hspace{7mm}=\frax{1}{(n+2^{n-2})^2}\left[n(n+2^{n-2})+(n+2)2^{n-2}\right]=$$
$$\hspace{-37,5mm}=\frax{n^2+(n+1)2^{n-1}}{(n+2^{n-2})^2}\,,$$as desired.
\hfill\rule{1,5mm}{1,5mm}
\smallskip

\bigskip\noindent{\bf Corollary 3.3.2.} {\it $\dd\lim_{n\to\infty}csd(Q_{2^n})=0.$}
\bigskip

The same type of reasoning will be used to calculate
$csd(S_{2^n})$.

\bigskip\noindent{\bf Theorem 3.3.3.} {\it The cyclic subgroup commutativity degree
of the quasi-dihedral group $S_{2^n}$ is
$$csd(S_{2^n})=\frax{n^2+3n\cdot2^{n-2}+5\cdot2^{n-3}}{(n+3\cdot2^{n-3})^2}\,.$$In particular, we have
$$csd(S_{16})=\frax{37}{50}\,.$$}
\smallskip

\noindent{\bf Proof.} It is a simple exercise to check that the
poset $L_1(S_{2^n})$ of cyclic subgroups of $S_{2^n}$ consists of
$$L(\langle x\rangle)\cup\{\langle x^{2k}y\rangle\mid k=0,1,...,2^{n-2}-1\}\cup\{\langle x^{2k+1}y\rangle\mid
k=0,1,...,2^{n-3}-1\}.$$Again, we have
$$|C_1(H)|=|L_1(S_{2^n})|=n+3\cdot2^{n-3},
\forall\hspace{1mm}H\leq\langle x\rangle.$$In order to study the
commutativity of the other two types of subgroups of $S_{2^n}$,
the following remarks are essential:
\begin{itemize}
\item[--] $\langle x^{2k_1}y\rangle\langle x^{2k_2}y\rangle=\langle x^{2k_2}y\rangle\langle x^{2k_1}y\rangle{\Longleftrightarrow} 2^{n-3}\mid k_1-k_2$;
\item[--] $\langle x^{2k_1+1}y\rangle\langle x^{2k_2+1}y\rangle=\langle x^{2k_2+1}y\rangle\langle x^{2k_1+1}y\rangle{\Longleftrightarrow} 2^{n-3}\mid k_1-k_2{\Longleftrightarrow} k_1=k_2$;
\item[--] $\langle x^{2k_1}y\rangle\langle x^{2k_2+1}y\rangle\neq\langle x^{2k_2+1}y\rangle\langle x^{2k_1}y\rangle, \forall\hspace{1mm} k_1,k_2$.
\end{itemize}We infer that
$$|C_1(\langle x^{2k}y\rangle)|=n+2, \forall\hspace{1mm}k=0,1,...,2^{n-2}-1$$and
$$|C_1(\langle x^{2k+1}y\rangle)|=n+1, \forall\hspace{1mm}k=0,1,...,2^{n-3}-1.$$Hence
\bigskip
$$\hspace{-42mm}csd(S_{2^n})=\frax{1}{(n+3\cdot2^{n-3})^2}\dd\sum_{H\in L_1(S_{2^n})}|C_1(H)|=$$
$$=\frax{1}{(n+3\cdot2^{n-3})^2}\left[\dd\sum_{^{H\in
L_1(S_{2^n})}_{\hspace{3mm}H\leq\langle
x\rangle}}\hspace{-3mm}|C_1(H)|+\hspace{-3mm}\dd\sum_{k=0}^{2^{n-2}-1}|C_1(\langle
x^{2k}y\rangle)|+\hspace{-3mm}\dd\sum_{k=0}^{2^{n-3}-1}|C_1(\langle
x^{2k+1}y\rangle)|\right]{=}$$
$$\hspace{22mm}=\frax{1}{(n+3\cdot2^{n-3})^2}\left[n(n+3\cdot2^{n-3})+(n+2)2^{n-2}+(n+1)2^{n-3}\right]=$$
$$\hspace{-43mm}=\frax{n^2+3n\cdot2^{n-2}+5\cdot2^{n-3}}{(n+3\cdot2^{n-3})^2}\,,$$completing
the proof.
\hfill\rule{1,5mm}{1,5mm}

\bigskip\noindent{\bf Corollary 3.3.4.} {\it $\dd\lim_{n\to\infty}csd(S_{2^n})=0.$}

\section{A criterion for a finite group to be Iwasawa}

A famous result by Gustafson \cite{3} concerning the commutativity degree states that if $d(G)>5/8$ then $G$ is abelian,
and we have $d(G)=5/8$ if and only if $G/Z(G)\cong\mathbb{Z}_2\times\mathbb{Z}_2$. In this section a similar problem
is studied for the cyclic subgroup commutativity degree, namely: \textit{is there a constant $c\in(0,1)$ such that if $csd(G)>c$ then
$G$ is Iwasawa}?
\bigskip

The answer to this problem is negative, as shows the following theorem.

\bigskip\noindent{\bf Theorem 4.1.} {\it The cyclic subgroup commutativity degree
of the non-Iwasawa group $\mathbb{Z}_{2^n}\times Q_8$, $n\geq 2$, tends to $1$ when $n$ tends to infinity.}
\bigskip

\noindent{\bf Proof.} Let $a=(a_1,a_2)\in \mathbb{Z}_{2^n}\times Q_8$. Then $o(a)=2^k$ if and only if either $o(a_1)=2^k$ and $o(a_2)\leq 2^k$ or $o(a_1)<2^k$ and $o(a_2)=2^k$.
We infer that $\mathbb{Z}_{2^n}\times Q_8$ has one element of order $1$, $3$ elements of order $2$, $28$ elements of order $4$, and $2^{k+2}$ elements of order $2^k$, $\forall\, k=3,4,...,n$.
These generate one cyclic subgroup of order $1$, $3$ cyclic subgroups of order $2$, $14$ cyclic subgroups of order $4$, and $8$ cyclic subgroups of order $2^k$, $\forall\, k=3,4,...,n$.
Consequently, $$|L_1(\mathbb{Z}_{2^n}\times Q_8)|=1+3+14+8(n-2)=8n+2.$$Then $$\barr{lcl} csd(\mathbb{Z}_{2^n}\times Q_8)&=&\frax1{(8n+2)^2}\,\left|\{(H,K)\in L_1(\mathbb{Z}_{2^n}\times Q_8)^2\mid
HK=KH\}\right|=\vsp &=&1-\frax1{(8n+2)^2}\,\left|\{(H,K)\in
L_1(\mathbb{Z}_{2^n}\times Q_8)^2\mid HK\neq KH\}\right|.\earr$$One the other hand, by Theorem 2.15 of \cite{1} we know that $\mathbb{Z}_{2^n}\times Q_8$ has $24(n+2)$ pairs of subgroups which do not permute. This implies that $$csd(\mathbb{Z}_{2^n}\times Q_8)\geq 1-\frax{24(n+2)}{(8n+2)^2}$$and so $\dd\lim_{n\to\infty}csd(\mathbb{Z}_{2^n}\times Q_8)=1$, completing the proof.
\hfill\rule{1,5mm}{1,5mm}

\bigskip\noindent{\bf Corollary 4.2.} {\it There is no constant $c\in(0,1)$ such that if $csd(G)>c$ then
$G$ is Iwasawa.}
\bigskip

However, we can get a positive answer to the above problem if we replace the condition "$csd(G)>c$" by the stronger condition "$csd^*(G)>c$",
where $$csd^*(G)={\rm min}\{csd(S) \mid S \mbox{ section of } G\}.$$This was suggested by the fact that a $p$-group is modular if and only
if each of its sections of order $p^3$ does. Moreover, if a $p$-group is not modular then it contains a section isomorphic to $D_8$ or to $E(p^3)$,
the non-abelian group of order $p^3$ and exponent $p$ for $p>2$ (see Lemma 2.3.3 of \cite{13}).

\bigskip\noindent{\bf Lemma 4.3.} {\it Let $G$ be a finite $p$-group such that $csd^*(G)>41/49$. Then $G$ is modular, and consequently an Iwasawa group.}
\bigskip

\noindent{\bf Proof.} Assume that $G$ is not modular. Then there is a section $S$ of $G$ such that $S\cong D_8$ or $S\cong E(p^3)$ for $p>2$. We can easily check that
$$csd(E(p^3))=\frax{p^3+5p^2+4p+4}{(p^2+p+2)^2}<\frax{41}{49}=csd(D_8).$$Therefore $csd(S)\leq 41/49$, contradicting our assumption.
\hfill\rule{1,5mm}{1,5mm}

\bigskip\noindent{\bf Lemma 4.4.} {\it Let $G$ be a finite group such that $csd^*(G)>19/25$. Then $G$ is nilpotent.}
\bigskip

\noindent{\bf Proof.} We will show by induction on $|G|$ that if $G$ is not nilpotent then $csd^*(G)\leq 19/25$, i.e. there is a section $S$ of $G$ with $csd(S)\leq 19/25$.
For $|G|=6$ we have $G\cong S_3$ and the desired conclusion follows by taking $S=G$. Assume now that it is true for all non-nilpotent
groups of order $<|G|$. We distinguish the following two cases.

If $G$ contains a proper non-nilpotent subgroup $H$, then $H$ has a section $S$ with $csd(S)\leq 19/25$ by the inductive hypothesis and we are done since
$S$ is also a section of $G$.

If all proper subgroups of $G$ are nilpotent, then $G$ is a Schmidt group. By \cite{12} (see also \cite{10}) it follows that $G$ is a solvable group of order $p^mq^n$
(where $p$ and $q$ are different primes) with a unique Sylow $p$-subgroup $P$ and a cyclic Sylow $q$-subgroup $Q$, and hence $G$ is a semidirect product of $P$ by $Q$.
Moreover, we have:
\begin{itemize}
\item[-] if $Q=\langle y\rangle$ then $y^q\in Z(G)$;
\item[-] $Z(G)=\Phi(G)=\Phi(P)\times\langle y^q\rangle$, $G'=P$, $P'=(G')'=\Phi(P)$;
\item[-] $|P/P'|=p^r$, where $r$ is the order of $p$ modulo $q$;
\item[-] if $P$ is abelian, then $P$ is an elementary abelian $p$-group of order $p^r$ and $P$ is a minimal normal subgroup of $G$;
\item[-] if $P$ is non-abelian, then $Z(P)=P'=\Phi(P)$ and $|P/Z(P)|=p^r$.
\end{itemize}We infer that $S=G/Z(G)$ is also a Schmidt group of order $p^rq$ which can be written as a semidirect product of an elementary abelian $p$-group $P_1$ of order $p^r$ by a cyclic group $Q_1$ of order $q$ (note that $S_3$ and $A_4$ are examples of such groups). Then $L_1(S)=L_1(P_1)\cup \{Q_1^x \mid x\in S\}$ and $$|L_1(S)|=\frax{p^r-1}{p-1}+1+p^r=\frax{p^{r+1}+p-2}{p-1}\,.$$One obtains: $$csd(S)=\frax{5p+4}{(p+2)^2}\hspace{3mm} \mbox{ for } r=1\0(a)$$and $$csd(S)=\frax{p^{2r}+3p^{r+2}-4p^{r+1}-p^r+p^2-4p+4}{(p^{r+1}+p-2)^2}\hspace{3mm} \mbox{ for } r\geq 2.\0(b)$$In both cases ($a$) and ($b$) we can easily check that $$csd(S)\leq\frax{19}{25}\,,$$as desired.
\hfill\rule{1,5mm}{1,5mm}
\bigskip

We are now able to prove the main result of this section.

\bigskip\noindent{\bf Theorem 4.5.} {\it Let $G$ be a finite group such that $csd^*(G)>41/49$. Then $G$ is an Iwasawa group. Moreover, we have
$csd^*(G)=41/49$ if and only if $G\cong G'\times G''$, where $G'$ is a $2$-group with $csd^*(G')=41/49$ and $G''$ is an Iwasawa group of odd order.}
\bigskip

\noindent{\bf Proof.} Since $csd^*(G)>41/49>19/25$, Lemma 4.4 implies that $G$ is nilpotent. Then it can be written as
$$G=\xmare{i=1}k G_i\,,\0(4)$$where $G_i$ is a Sylow $p_i$-subgroup of $G$, $i=1,2,...,k$. For each $i$ we have $$csd^*(G_i)\geq csd^*(G)>\frax{41}{49}\,,$$and therefore
$G_i$ is an Iwasawa group by Lemma 4.3. Consequently, $G$ is also an Iwasawa group.

Suppose now that $csd^*(G)=41/49$. Then $G$ is nilpotent by Lemma 4.4, and therefore it has a direct decomposition of type $(4)$, where we can assume $p_1<p_2<...<p_k$. Remark that
$p_1=2$. Indeed, if $p_1>2$ then all $p_i$'s are odd, which implies that $G_i$ cannot have sections isomorphic with $D_8$, $\forall\, i=1,2,...,k$. On the other hand, $G_i$ cannot also have sections isomorphic with $E(p_i^3)$ because $csd^*(G_i)\geq csd^*(G)=41/49$. Thus $G_i$ is Iwasawa, $\forall\, i=1,2,...,k$, and the same thing can be said about $G$, a contradiction. Hence $p_1=2$ and we are done by taking $$G'=G_1 \mbox{ and }  G''=\xmare{i=2}k G_i.$$

Conversely, since $G'$ and $G''$ are of coprime orders, every section $S$ of $G\cong G'\times G''$ is of type $S\cong S'\times S''$, where $S'$ and $S''$ are sections of $G'$ and $G''$, respectively. Then $$csd(S)=csd(S')csd(S'')=csd(S')$$because $G''$ is Iwasawa. This shows that $$csd^*(G)=csd^*(G')=\frax{41}{49}\,,$$completing the proof.
\hfill\rule{1,5mm}{1,5mm}
\bigskip

We end this section by noting that the problem of finding the structure of $2$-groups $G'$ with $csd^*(G')=41/49$ remains open.

\section{Conclusions and further research}

Similarly with our previous concepts of \textit{subgroup commutativity
degree}, \textit{normality degree} or \textit{cyclicity degree} of
a finite group, the \textit{cyclic subgroup commutativity degree}
can also constitute a significant aspect of probabilistic finite
group theory. Clearly, the study started in this paper can
successfully be extended to other classes of finite groups and all
problems on $sd(G)$, $ndeg(G)$, $cdeg(G)$ (see e.g. \cite{16}-\cite{21})
can be investigated for $csd(G)$, too. On the other hand, the connections
between the above concepts seem to be very inte\-resting. These will surely
constitute the subject of some further research.
\smallskip

Finally, we formulate several specific open problems on cyclic
subgroup commutativity degrees.

\bigskip\noindent{\bf Problem 5.1.} Compute explicitly the cyclic
subgroup commutativity degree of ${\rm ZM}(m,n,r)$ (see
\cite{18}), or, more generally, the cyclic subgroup
commutativity degree of an arbitrary metacyclic group.

\bigskip\noindent{\bf Problem 5.2.} Let $G$ be a finite group. Study the properties
of the map $csd: L(G)\longrightarrow [0,1]$, $H\mapsto csd(H)$. Is
it true that for every $H,K\in L(G)$, we have $H\subseteq
K\Longrightarrow csd(H)\geq csd(K)$?

\bigskip\noindent{\bf Problem 5.3.} For many finite groups $G$, the commutativity
of $x,y\in G$ is strongly connected with the commutativity of
$\langle x\rangle, \langle y\rangle\in L_1(G)$. Can be extended
this to a connection between $d(G)$ and $csd(G)$?

\bigskip\noindent{\bf Problem 5.4.} Does exist finite groups $G$
such that $csd(G)=sd(G)\neq 1$?

\vspace*{3ex}
\small

\begin{minipage}[t]{7cm}
Marius T\u arn\u auceanu \\
Faculty of  Mathematics \\
"Al.I. Cuza" University \\
Ia\c si, Romania \\
e-mail: {\tt tarnauc@uaic.ro}
\end{minipage}
\hfill
\begin{minipage}[t]{7cm}
Mihai-Silviu Lazorec \\
Faculty of  Mathematics \\
"Al.I. Cuza" University \\
Ia\c si, Romania \\
e-mail: {\tt Mihai.Lazorec@math.uaic.ro}
\end{minipage}
\end{document}